\documentclass{amsart}[12pt]
\def\doctype{}

\usepackage{latexsym,amssymb,amsfonts,dsfont,bm}
\usepackage{fancyhdr,caption}
\usepackage{tikz}
\usepackage{hyperref}

\newcommand\lam{\lambda}

\newcommand{\cB}{\mathcal{B}}

\newcommand{\comment}[1]{}

\numberwithin{equation}{section}

\fancypagestyle{titlepage}{
\fancyhead[R]{\doctype}
\fancyhead[CL]{}
\cfoot{\vspace{5pt} \thepage}
}

\let\oldsection\section
\newcommand\boldsection[1]{\oldsection{\bf #1}}
\newcommand\starsection[1]{\oldsection*{\bf #1}}
\makeatletter
\renewcommand\section{\@ifstar\starsection\boldsection}
\makeatother

\newtheoremstyle{theorem}
  {12pt}		  
  {0pt}  
  {\sl}  
  {\parindent}     
  {\bf}  
  {. }    
  { }    
  {}     
\theoremstyle{theorem}
\newtheorem{thm}{Theorem}[section]  
\newtheorem{lemma}[thm]{Lemma}     

\newtheorem{cons}[thm]{Construction}
\newtheorem{prop}[thm]{Proposition}

\newtheoremstyle{definition}
  {12pt}		  
  {0pt}  
  {}  
  {\parindent}     
  {\bf}  
  {. }    
  { }    
  {}     
\theoremstyle{definition}

\newtheorem{ex}[thm]{Example}

\newcommand\rks{{\sc Remarks.} }

\renewcommand{\proofname}{Proof}

\makeatletter
\renewenvironment{proof}[1][\proofname]{\par
  \pushQED{\qed}%
  \normalfont \partopsep=\z@skip \topsep=\z@skip
  \trivlist
  \item[\hskip\labelsep
        \scshape
    #1\@addpunct{.}]\ignorespaces
}{%
  \popQED\endtrivlist\@endpefalse
}
\makeatother

\makeatletter
\renewcommand*\@maketitle{%
  \normalfont\normalsize
  \@adminfootnotes
  \@mkboth{\@nx\shortauthors}{\@nx\shorttitle}%
  \global\topskip42\p@\relax 
  \@settitle
  \ifx\@empty\authors \else {\vskip 1em
\vtop{\centering\shortauthors\@@par}} \fi
  \ifx\@empty\@date \else {\vskip 1em \vtop{\centering\@date\@@par}}\fi 
  \ifx\@empty\@dedicatory
  \else
    \baselineskip18\p@
    \vtop{\centering{\footnotesize\itshape\@dedicatory\@@par}%
      \global\dimen@i\prevdepth}\prevdepth\dimen@i
  \fi
  \@setabstract
  \normalsize
  \if@titlepage
    \newpage
  \else
    \dimen@34\p@ \advance\dimen@-\baselineskip
    \vskip\dimen@\relax
  \fi
} 
\renewcommand*\@adminfootnotes{%
  \let\@makefnmark\relax  \let\@thefnmark\relax
  \ifx\@empty\@subjclass\else \@footnotetext{\@setsubjclass}\fi
  \ifx\@empty\@keywords\else \@footnotetext{\@setkeywords}\fi
  \ifx\@empty\thankses\else \@footnotetext{%
    \def\par{\let\par\@par}\@setthanks}%
  \fi
\thispagestyle{titlepage}
}
\makeatother

\title{\large Constructions of Sarvate-Beam Group Divisible Designs}

\author{Peter J.~Dukes}
\address{\rm Mathematics and Statistics, 
University of Victoria, Victoria, BC, Canada}
\email{dukes@uvic.ca}

\author{Joanna Niezen}
\address{\rm Mathematics and Statistics, 
Simon Fraser University, Burnaby, BC, Canada}
\email{jniezen@sfu.ca}

\thanks{Research of Peter Dukes is supported by NSERC grant 312595--2017}

\date{\today}

\begin{document}

\begin{abstract}
A balanced incomplete block design is a set system in which all pairs of distinct elements occur with a constant frequency.  By contrast, a Sarvate-Beam design induces an interval of distinct frequencies on pairs.  In this paper, we settle the existence of a Sarvate-Beam variant of group divisible designs of uniform type with block size three.
\end{abstract}

\maketitle
\hrule


\section{Introduction}

Given a set system $(V,\cB)$,  let us say that the \emph{frequency} of $X \subseteq V$ equals 
$\sum_{B \in \cB} \mathds{1}_{X \subseteq B}$, the number of members of $\cB$ (counting multiplicity) that contain $X$.  A $t$-design is a set system in which all sets of size $i \le t$ have a constant frequency $\lam_i$.  An important special case is that of a \emph{Steiner triple system}, where $|B|=3$ for every $B \in \cB$ and all two-element subsets have frequency 1.  Letting $v=|V|$, a Steiner triple system on $v$ elements is known to exist \cite{Kirkman} if and only if $v \equiv 1,3 \pmod{6}$.  As in this case, our attention in what follows restricted to frequencies of 2-subsets.

For the recursive construction of designs or for certain experimental applications, it can be useful to have some  frequencies equal to zero.  With this in mind, a {\emph{group divisible design} is a triple $(V,\Pi,\cB)$, where $V$ is a set equipped with a partition $\Pi=\{V_1,\dots,V_u\}$ into `groups', and 
$\cB$ is a collection of subsets of $V$ called `blocks', such that 

\vspace{-5mm}
\begin{itemize}
\item
every pair of elements in the same group has frequency zero, and
\item
every pair of elements in different groups has frequency one.
\end{itemize}
This is abbreviated as a $k$-GDD of type $g_1^{u_1}g_2^{u_2}\cdots g_s^{u_s}$ when the block size is a constant $k$ and there are exactly $g_i$ groups of size $u_i$, $i=1,\dots,s$.  Here, we are normally interested in the special case of `uniform' group size; that is, GDDs of type $g^u$.

Sarvate and Beam \cite{SB1} introduced a new condition for set systems in which all 2-subsets have a {\em different} frequency. They called such a set system an `adesign', focusing the definition on constant block size $k$.  With no additional restrictions, it is usually easy to construct adesigns.  A more challenging special case, now called a {\em Sarvate-Beam design} SB$_{\mu}(v,k)$, demands that the pair frequencies form an interval $\{\mu,\mu+1,\dots,\mu+\binom{v}{2}-1\}$ of distinct values. Since the total of the pairwise frequencies is equal to $\binom{k}{2}$ times the number of blocks, a necessary divisibility for the existence of an SB$_\mu(v,k)$ is

\vspace{-5mm}
\begin{equation}
\label{sb-div}
\tbinom{k}{2} \text{~~divides~~} \mu \tbinom{v}{2}+\tfrac{1}{2} \tbinom{v}{2} \left( \tbinom{v}{2}-1 \right).
\end{equation}

An SB$_{\mu}(v,3)$ is called a {\em Sarvate-Beam triple system} and abbreviated SBTS$_{\mu}(v)$. With $k=3$, condition \eqref{sb-div} reduces to 
$v \equiv 0,1 \pmod{3}$ or $\mu \equiv 0 \pmod{3}$.
Ma, Chang and Feng in \cite{MCF} settled the existence of Sarvate-Beam triple systems in the case $\mu=1$ for all admissible $v>3$, although they used the term `strictly pairwise distinct triple system' and the notation sPDTS$(v)$ to mean what we call an SBTS$_1(v)$.  Then, Dukes and Short-Gershman \cite{DSG} proved that the necessary conditions for existence of SBTS$_\mu(v)$ are sufficient for arbitrary $\mu \ge 0$, with the exception $(v,\mu)=(4,0)$.


Our focus in this paper is a generalization inspired by GDDs.
A \emph{Sarvate-Beam group divisible design} (SBGDD) is a triple $(V,\Pi,\cB)$ where $V$ is a set, $\Pi=\{V_1,\dots,V_u\}$ is a partition of $V$ into `groups' and $\cB$ is a multiset of subsets of $V$ called `blocks', such that

\vspace{-5mm}
\begin{itemize}
\item
every pair of elements in the same group has frequency zero, and
\item
the frequencies of $2$-element subsets across different groups are consecutive integers.
\end{itemize}

We restrict our attention to the case in which all blocks have size three.
The \emph{type} of an SBGDD is the integer partition of $v=|V|$ induced by the group sizes $|V_i|$.  Of particular interest is the uniform case where all group sizes are equal, say $|V_i|=g$ for each $i$.   Following the convention used in (ordinary) GDDs, we use exponential notation $g^u$ to abbreviate the type of an SBGDD.  To emphasize that the starting frequency equals $\mu$, we use the notation SBGDD$_\mu$. 

In slightly different language, existence of an SBGDD$_\mu$ of type $g^u$ is equivalent to a multiset of triangles in the complete multipartite graph $K_u \cdot \overline{K_g}$ covering the edges with frequencies $\mu,\mu+1,\mu+2,\dots,\mu+g^2 \binom{u}{2}-1$. 
The existence of SBGDDs was first considered by Moolsombut and Hemakul \cite{MH}, where a few small examples were given.

We begin with the important case $\mu=0$.  Note that the total of all pair frequencies equals 
$$\frac{1}{2} g^2 \tbinom{u}{2} \left(g^2 \tbinom{u}{2}-1 \right) \equiv 0 \pmod{3}$$
for any integers $g$ and $u$ since both $g^2$ and $\binom{u}{2}$ are either $0$ or $1 \pmod{3}$.

From the definition, $u \ge 3$ is necessary for the existence of an SBGDD of type $g^u$.  The boundary case $u=3$ has an interesting alternate interpretation.
A \emph{Sarvate-Beam cube} SBC$(n)$ is an $n \times n \times n$ array in which each cell contains a nonnegative integer and such that the $3n^2$ line sums cover the interval of integers $\{0,1,\dots,3n^2-1\}$.
The above definition was introduced in \cite{DSG}, and example cubes were given for $n=2,3$.  In \cite{DN}, it was shown that SBC$(n)$ exist for every positive integer $n \ge 2$.  An SBC$(n)$ is equivalent to a SBGDD of type $n^3$ in the following way: the $(a,b,c)$-entry of the SBC gives the frequency of triple $\{(a,1),(b,2),(c,3)\}$ on the point set $V=[n] \times [3]$, where group partition $\Pi = \{[n] \times \{i\}: i = 1,2,3\}$ is used.  It is clear that the line sums through the cube correspond to the total frequency of pairs from different groups.

Our main result decides the existence of uniform SBGDDs with starting frequency zero.

\begin{thm}
\label{main}
For integers $g \ge 1$ and $u \ge 3$, there exists an SBGDD$_0$ of type $g^u$, with the exception of $(g,u)=(1,4)$.
\end{thm}

As mentioned above, the case $g=1$ was previously settled in \cite{DSG} and the case 
$u=3$ was recently settled in \cite{DN}.  The result here is a common generalization and the proof is constructive.

The rest of the paper is arranged as follows.  In the next section, we review some background on group divisible designs, especially those with mixed block sizes.  Then, in Section~\ref{sec:recursive}, we give a broad recursive construction for SBGDDs that carries most of the subsequent results.  The idea behind this construction is actually very simple: we obtain an SBGDD if its underlying graph can be edge-decomposed into tiles, each of which realizes a sub-interval of the needed range of frequencies.
In Section~\ref{sec:special}, some special `atomic' constructions are given.  For the proof of the main result, which occurs in Section~\ref{sec:combining}, these seed designs are inflated and/or combined to tackle the cases $u=4,5,6,8$, at which point our recursion is relatively easily able to settle the remaining values of $u$.  Finally, in Section~\ref{sec:mu}, we give some minor modifications of our constructions to handle larger starting frequencies $\mu > 0$.


\section{Background on group divisible designs}

This section collects some background existence results on group divisible designs that will come in handy later.
We begin by discussing the extremal case in which every block traverses the group partition.  A \emph{transversal design} TD$(k,n)$ is a $k$-GDD of type $n^k$.  A TD$(k,n)$ is equivalent to a set of $k-2$ mutually orthogonal latin squares (MOLS) of order $n$; see for instance \cite[III.3.18]{handbook}.  In particular, there exists a finite field construction when $n \ge k-1$ is a prime power.  From this and the data on MOLS in \cite[III.3.81]{handbook} we obtain existence in many cases when $k$ is small.  The following is enough for our purposes.

\begin{lemma}[See \S III of \cite{handbook}]
\label{td-existence}~

\vspace{-5mm}
\begin{enumerate}
\item
There exists a TD$(k,n)$ if $n$ has a factorization into prime powers, each at least $k-1$.
\item
For $k \in \{3,4,5\}$, there exists a TD$(k,n)$ provided $(k,n) \not\in \{(4,2),(4,6),(5,2),(5,6),(5,10)\}$.
\end{enumerate}
\end{lemma}

Next, it is useful to consider group divisible designs with mixed block sizes.
If the set of allowed block sizes in a GDD is $K \subseteq \{2,3,\cdots\}$, the abbreviation $K$-GDD is used.  Note that it is not necessary for every integer in $K$ to be realized as a block size, but only that block sizes belong to $K$.  In the case that $|V|=v$ and $\Pi$ is the partition into singletons, we obtain a \emph{pairwise balanced design}, or PBD$(v,K)$.  Such designs are foundational for R.M.~Wilson's classical existence theory \cite{W1,W2} for $2$-designs.

\begin{lemma}
\label{small-gdds}
The following GDDs exist:

\vspace{-3mm}
\begin{tabular}{ll}
{\rm (a)} $3$-GDD of type $2^4$
&
{\rm (e)} $3$-GDD of type $2^6$
\\
{\rm (b)} 
$\{3,5\}$-GDD of type $2^5$
&
{\rm (f)} $\{3,5\}$-GDD of type $2^8$
\\
{\rm (c)} $4$-GDD of type $3^5$
&
{\rm (g)} $4$-GDD of type $3^8$
\\
{\rm (d)} $\{3,4\}$-GDD of type $3^6$
&
{\rm (h)} $\{4,5\}$-GDD of type $5^8$.
\end{tabular}
\end{lemma}

\begin{proof}
The cases with a single block size are straightforward point-deletions of finite geometries or other known small designs.  We highlight the more subtle constructions.

(b) Start from an incomplete triple system of order 11 with a hole $H$ of order 5.  (We refer the reader to \cite{handbook,ts} for a formal definition and constructions.)  Delete an element $x \not\in H$, and remove all blocks incident with $x$. The removal of these blocks induces a partition of the remaining points into 5 groups, each of size 2.  All other blocks are kept, and $H$ is turned into a block.

(d) From a Kirkman triple system of order $15$ (see again \cite{handbook,ts}), turn one parallel class into a set of groups, producing a resolvable $3$-GDD of type $3^5$.  Now, extend three parallel classes into blocks of size four, each joined to a new point in a sixth group.

(f) This is similar to (b), except starting from an incomplete triple system of order 17.

(h) Delete one group $G$ from a $5$-GDD of type $5^9$, one example of which appears in \cite[VI.16.31]{handbook}.  Any blocks which intersect $G$ are reduced to have size four. 
\end{proof}

It is useful to have various group divisible designs that facilitate `lifting' an interval of frequences to achieve larger starting values.  Let $\lam$ be a nonnegative integer.  A GDD of \emph{index} $\lambda$ is a triple $(V,\Pi,\cB)$, where $V$ and $\Pi$ are as before, and $\cB$ is a collection of blocks such that

\vspace{-5mm}
\begin{itemize}
\item
every pair of elements in the same group has frequency zero, and
\item
every pair of elements in different groups has frequency $\lam$.
\end{itemize}

If all blocks have size three, this is abbreviated as a $(3,\lam)$-GDD.  The group type is abbreviated with the usual exponential notation.

\begin{ex}
An SBGDD$_\mu$ of type $2^3$ can be constructed as the multiset union of blocks of an SBGDD$_0$ of type $2^3$ (equivalent to a Sarvate-Beam cube of order $2$) and a $(3,\mu)$-GDD on the same points and groups (obtainable from $\mu$ copies of a $3$-GDD of type $2^3$).
\end{ex}

It is easy to see that $3$-GDDs of type $m^3$ exist for every positive integer $m$, as these are equivalent to $m \times m$ latin squares.  The next result gathers some results for $3$-GDDs having four, five, six, or eight groups.    The object referenced in part (c) is called an `incomplete group divisible design'.  For our purposes it is enough to say that this is equivalent to a $\{3,5\}$-GDD of type $m^5$ with a single block of size $5$ that is removed from the block family $\cB$.
The reader is referred to \cite[IV.1.37]{handbook} for more details on IGDDs.

\begin{lemma}[See \cite{Zhu}]\label{zhu}
Let $m$ and $\lambda$ be positive integers.  There exists:

\vspace{-5mm}
\begin{enumerate}
\item
a $(3,\lambda)$-GDD of type $m^4$ if $m$ or $\lambda$ is even;
\item
a $(3,\lambda)$-GDD of type $m^5$ if $3 \mid m$ or $3 \mid \lambda$;
\item
a $3$-IGDD of type $(m;1)^5$ if $3 \nmid m$.

\item
a $(3,\lambda)$-GDD of type $m^6$ if $m$ or $\lambda$ is even;
\item
a $(3,\lambda)$-GDD of type $m^8$ if $m$ is even, and $3 \mid m$ or $3 \mid \lambda$.
\end{enumerate}
\end{lemma}

\section{Recursive constructions}
\label{sec:recursive}

The recursive constructions for Sarvate-Beam GDDs to follow are similar in spirit to similar constructions for classical block designs.  The key difference is that the intervals of pair frequencies in ingredient SBGDDs must be abutted to form a longer interval in the resultant SBGDD.  We begin with an illustrative example.

\begin{ex}[SBGDD of type $2^7$]
\label{sbgdd2-7}
Consider the Fano plane on elements $V=\{0,1,\dots,6\}$ whose blocks are the (mod 7) translates of $\{0,1,3\}$.  We construct an SBGDD of type $2^7$ on the set $V \times \{1,2\}$.  Each element $x$ is replaced by two elements $(x,1)$ and $(x,2)$.  Every block of the design is replaced by an SBGDD of type $2^3$.  The list of frequencies so produced is $0^7,1^7,\dots,11^7$, where there are seven occurrences of each frequency in the ingredient design.  In order to have the frequency lists form an interval, we `lift' the frequencies from the $i$th block $\{i,1+i,3+i\}$ by $12i$, $i=0,1,\dots,6$.  This can be done by appending $12i$ copies of a GDD of type $2^3$ (equivalent to a latin square of side $2$) on the corresponding points.  After doing so, on $\{i,1+i,3+i\} \times \{1,2\}$, pairs of elements (with distinct first coordinates) occur in the resultant design with frequencies $12i,12i+1,\dots,12i+11$.  In this way, the full list of frequencies obtained is $0,1,\dots,12\times 6+11$, as desired.
\end{ex}

The construction of Example~\ref{sbgdd2-7} is reminiscent of Wilson's fundamental construction, \cite{W1}, where small ingredient GDDs are aligned on blocks of a master design to produce a larger GDD.  A technicality for Sarvate-Beam GDDs, however, is that the ingredients must have the correct starting frequencies, either by lifting via an (ordinary) GDD, as in the example, or otherwise via a direct construction.

There are two situations we encounter in the sequel where this technicality requires careful attention.  The first occurs when we wish to replace blocks of size four with SBTS$_\mu(4)$.  If every block of the master design has size four, there is no way to choose a starting block to obtain frequencies $0,1,\dots,5$, since an SBTS$_0(4)$ does not exist.  In this situation, we make a small adjustment to the ordering of frequencies at the beginning.  Choose five blocks of the form $\{a,b,c,d\}$, $\{a,x,*,*\}$, $\{b,x,*,*\}$, $\{c,y,*,*\}$, $\{d,y,*,*\}$, where it is not important that the starred points be distinct across the different blocks.  On $\{a,b,c,d\}$, instead of an SBTS$(4)$, we use the multiset of triples $\{a,b,c\}$, $\{b,c,d\}$, $\{b,c,d\}$, which induces pair frequencies $0,1,1,2,2,3$.
Next, include three copies of $\{a,b,x\}$ and three copies of $\{c,d,y\}$, so that $\{a,b\}$ has frequency $4$ and $\{c,d\}$ has frequency $5$.  This gives the correct interval of frequencies on $\{a,b,c,d\}$.  Now, the other four blocks are assigned triples so as to achieve the remaining frequencies from $0$ to $29$.  An illustration is shown in Figure~\ref{starting4}.  After handling these first five blocks, the remaining blocks can be replaced as usual by SBTS$_\mu(4)$ to produce the desired SBGDD.

\begin{figure}[htbp]
\begin{center}
\begin{tikzpicture}[scale=2]
\draw[blue] (0,0)--(1,0);
\draw[blue] (0,0.04)--(1,0.04);
\draw[red] (0,-0.04)--(1,-0.04);
\draw[red] (0,0)--(0,1);
\draw[red] (1,0)--(0,1);
\draw[blue] (0.015,-0.015)--(1.015,0.985);
\draw[blue] (-0.015,0.015)--(0.985,1.015);
\draw[blue] (.98,0)--(.98,1);
\draw[blue] (1.02,0)--(1.02,1);

\filldraw[black!10] (-0.03,0)--(-.8,.5)--(-0.03,1);
\filldraw[black!10] (1.03,0)--(1.8,.5)--(1.03,1);

\draw (0,0)--(-.8,.5)--(-1.3,-.3)--(-.5,-.8)--(0,0);
\draw (0,1)--(-.8,.5)--(-1.3,1.3)--(-.5,1.8)--(0,1);
\draw (1,0)--(1.8,.5)--(2.3,-.3)--(1.5,-.8)--(1,0);
\draw (1,1)--(1.8,.5)--(2.3,1.3)--(1.5,1.8)--(1,1);

\node at (0.5,1) {$0$};
\node at (0.3,0.7) {$1$};
\node at (0.3,0.3) {$2$};
\node at (0.5,-0.15) {$3$};

\node at (-.35,.5) {\Large $+3$};
\node at (1.25,.5) {\Large $+3$};

\node[align=center] at (-.65,-.15) {$5,6,7,$ \\ $9,10,11$};
\node[align=center] at (-.65,1.15) {$15,16,17,$ \\ $11,12,13$};
\node[align=center] at (1.65,-.15) {$23,24,25,$ \\ $27,28,29$};
\node[align=center] at (1.65,1.15) {$21,22,23,$ \\ $17,18,19$};

\draw[->] (-.81,.07)--(-.45,.25);
\draw[->] (-.85,.95)--(-.45,.75);
\draw[->] (1.4,.07)--(1.45,.25);
\draw[->] (1.4,.96)--(1.45,.75);

\filldraw (0,0) circle [radius=.05];
\filldraw (1,0) circle [radius=.05];
\filldraw (0,1) circle [radius=.05];
\filldraw (1,1) circle [radius=.05];
\filldraw (-.8,.5) circle [radius=.05];
\filldraw (-.5,-.8) circle [radius=.05];
\filldraw (-1.3,-.3) circle [radius=.05];
\filldraw (-.5,1.8) circle [radius=.05];
\filldraw (-1.3,1.3) circle [radius=.05];
\filldraw (1.8,.5) circle [radius=.05];
\filldraw (1.5,-.8) circle [radius=.05];
\filldraw (2.3,-.3) circle [radius=.05];
\filldraw (1.5,1.8) circle [radius=.05];
\filldraw (2.3,1.3) circle [radius=.05];
\end{tikzpicture}
\end{center}
\caption{Starting frequencies using blocks of size four}
\label{starting4}
\end{figure}
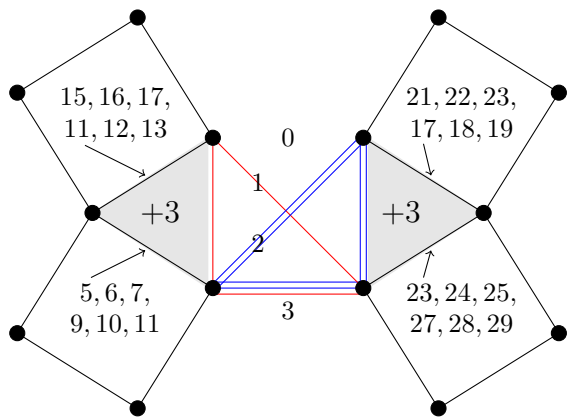

The other case which requires some attention is when we wish to replace more than one block of size $5$ with SBTS$_\mu(5)$, which do not exist unless $3 \mid \mu$.  An SBTS$(5)$ placed on the first such block block would cover frequencies $0,1,\dots,9$, but the next could not begin at $10$.  To overcome this problem, we make a small shuffle to the frequency lists covered by the next blocks.  Observe that the configuration of triples $\{a,b,c\}$, $\{a,d,e\}$, $\{b,d,e\}$, $\{c,d,e\}$ induces pair frequencies $1,1,\dots,1,3$.  And, the complete $(5,3,3)$-design on points $a,b,c,d,e$ lets us lift frequencies on these pairs by a multiple of three.  So, taking the multiset union of blocks of an SBTS$(5)$, three complete designs, and the configuration above, we can produce frequency list $10,11,\dots,18,21$.  Next, taking nine copies of the configuration in which exactly one pair is omitted as the threefold covered pair, we can produce frequency list $9,11,\dots,11$.  
Take the multiset union with SBTS$(5)$ so that the pair of frequency $1$ is aligned with the pair of frequency $9$.  This produces frequency list $10,11,13,14,\dots,20$.  Lifting as before with three complete designs, we obtain frequency list $19,20,22,23,\dots,29$.  Using this variant on lifting frequencies, we can take blocks of size five three at a time and form a long interval of frequencies.  Table~\ref{three5} summarizes the above combinations of frequencies.  A similar tactic can be applied to blocks of size 8; see \cite{Joey} for details.  In all applications to follow, the number of blocks of size $5$ or $8$ is congruent to $0$ or $1 \pmod{3}$; that is, we are never left with an unpatched hole in the frequency list.

\begin{table}[htbp]
$$\begin{array}{c|cccccccccc}
\text{pair} & ab & ac & ad & ae & bc & bd & be & cd & ce & de \\
\hline
\text{SBTS} & 0 & 1 & 3 & 2 & 5 & 4 & 7 & 6 & 8 & 9\\
\text{$(5,3,9)$-design} & 9 & 9 & 9 & 9 & 9 & 9 & 9 & 9 & 9 & 9\\
\text{one configuration} & 1 & 1 & 1 & 1 & 1 & 1 & 1 & 1 & 1 & 3\\
\hline
\text{totals} & 10 & 11 & 13 & 12 & 15 & 14 & 17 & 16 & 18 & 21\\
\end{array}$$

$$\begin{array}{c|cccccccccc}
\text{pair} & ab & ac & ad & ae & bc & bd & be & cd & ce & de \\
\hline
\text{SBTS} & 0 & 1 & 3 & 2 & 5 & 4 & 7 & 6 & 8 & 9\\
\text{$(5,3,9)$-design} & 9 & 9 & 9 & 9 & 9 & 9 & 9 & 9 & 9 & 9\\
\text{~9 configurations} & 11 & 9  & 11 & 11 & 11 & 11 & 11 & 11 & 11 & 11\\
\hline
\text{totals} & 20 & 19 & 23 & 22 & 25 & 24 & 27 & 26 & 28 & 29\\
\end{array}$$

\caption{Variant on lifting frequencies for blocks of size five}
\label{three5}
\end{table}

Our main recursive construction, which gives an analogue of Wilson's fundamental construction, can be stated in the language of graphs.  Roughly speaking, it says that an SBGDD can be built `by tiles' which decompose the edges of the underlying complete multipartite graph.

\begin{cons}
\label{main-cons}
Suppose the complete multipartite graph $K_u \cdot \overline{K_g}$ has an edge-decomposition into subgraphs $G_1,G_2,\dots,G_k$ and, for each $G_i$, there exists a multiset of triangles which cover the edges of $G_i$ with frequency list $[s_{i-1},s_i-1]$, where $s_0=\mu$ and in general $s_i = \mu+\sum_{j<i} |E(G_j)|$.  Then there exists an SBGDD$_\mu$ of type $g^u$.
\end{cons}

Often in our use of Construction~\ref{main-cons}, the graphs $G_i$ are simply taken as cliques, for instance from blocks in a master design, to be replaced in the construction by SBTS.  More generally, we can let $G_i$ be complete multipartite graphs and replace with known SBGDDs.  Referencing the special situations above, we can if necessary choose some $G_i$ as an edge-disjoint union of five cliques of size four (see Figure~\ref{starting4}), or three edge-disjoint cliques of order 5 or 8.

\section{Special constructions}
\label{sec:special}

This section contains some small examples of SBGDDs which are useful as ingredients to our recursive constructions.  The first few with group size $g=2$ were found with computer-assisted search.


\begin{ex}[SBGDD of type $2^4$; see also \cite{DSG}]
With group partition $\{\{a,a'\},\{b,b'\},\{c,c'\},\{d,d'\}\}$, the following block multiplicities form an SBGDD of type $2^4$.
\begin{center}
\begin{tabular}{cccccccc}
$ab'c$&	$ab'c'$&	$2ab'd$&	$2ab'd'$&	$3acd$&	$4acd'$&	$4ac'd$&	$5ac'd'$\\
$a'bd'$&	$a'b'd$&	$a'b'd'$&	$a'cd$&	$2a'cd'$&	$2a'c'd$&	$3a'c'd'$&	$6bcd$\\
$8bcd'$&	$6bc'd$&	$7bc'd'$&	$7b'cd$&	$7b'cd'$&	$10b'c'd$&	$8b'c'd'$&	
\end{tabular}
\end{center}
\end{ex}

\begin{ex}[SBGDD of type $2^5$ with a sub-SBTS$(5)$]
With group partition $\{\{a,a'\},\dots,\{e,e'\}\}$, the following block multiplicities result in frequencies $10,11,\dots,39$ on pairs of the form $\{x,y'\}$ or $\{x',y'\}$.
\begin{center}
\begin{tabular}{ccccccccc}
$10ae'c'$&	$3ae'b'$&	$ac'b'$&	$10ad'b'$&	$2be'c'$&	$11be'a'$&	$4be'd'$&	$7bc'a'$&	$15bc'd'$\\
$8ba'd'$&	$9ce'a'$&	$13ce'd'$&	$6ce'b'$&	$3ca'b'$&	$12cd'b'$&	$7de'c'$&	$9de'b'$&	$5dc'a'$\\
$7dc'b'$&	$13da'b'$&	$6ec'a'$&	$6ec'd'$&	$8ec'b'$&	$13ea'd'$&	$4ea'b'$&	$3ed'b'$&	$4e'c'a'$\\
$6e'c'd'$&	$6e'c'b'$&	$5e'a'd'$&	$5e'a'b'$&	$3e'd'b'$&	$6c'a'd'$&	$5c'a'b'$&	$3c'd'b'$&	$7a'd'b'$
\end{tabular}
\end{center}
If we take the multiset union with an SBTS$(5)$ on $\{a,b,\dots,e\}$, the result is an SBGDD of type $2^5$.  Moreover, using copies of a $3$-IGDD of type $(2;1)^5$, we can raise the starting frequency of our `incomplete' SBGDD to any integer $\mu \ge 10$.
\end{ex}

\begin{ex}[SBGDD of type $2^6$ with a sub-SBTS$(6)$]
Similar to the previous example, we can begin with an SBTS$(6)$ on $\{a,b,\dots,f\}$, and then realize larger frequencies starting at $\binom{6}{2}=15$ using blocks of the form $\{x,y',z'\}$ or $\{x',y',z'\}$.  There is considerable freedom in assigning block multiplicities to accomplish this: for instance, one can begin with $7ab'c'$, $8ab'd'$, $9ab'e'$ so that $ab'$ has frequency $7+8=15$, and $ac'$ has frequency $7+9=16$, and so on.  Multiplicities for blocks of the form $\{x',y',z'\}$ to finish off the construction are found without much difficulty by computer.
More details on the methodology can be found in the second author's dissertation \cite{Joey}. See the appendix for a complete list of block multiplicities; the first row of the block list gives a sub-SBTS$(6)$.
\end{ex}

\begin{ex}[SBGDD of type $2^8$]
This construction uses the existence of an SBGDD of type $2^6$ with two appended groups of size two. An SBC$(2)$ is placed on the two new groups with each of the original groups. Unfortunately the pair frequencies do not abut perfectly, as the points between the two added groups are covered in the blocks of every Sarvate-Beam cube. In order to reduce the four pair frequencies of the added points, ad-hoc block swaps are performed. The blocks and multiplicities are listed in the appendix. 
\end{ex}

We note a few other useful small SBGDDs, where now we are able to begin to use Construction~\ref{main-cons} to assemble them in pieces.

\begin{prop} \label{prop3^u}
There exists an SBGDD of type $3^u$ for each $u \in \{4,5,8\}$.
\end{prop}

\begin{proof}
Referring to Lemma~\ref{small-gdds}, begin with a $4$-GDD of type $3^u$.  Apply Construction~\ref{main-cons}, replacing blocks with SBTS$_\mu(4)$.  Here, we take care to reserve five blocks of size $4$ for starting the interval of frequencies, as described in Section~\ref{sec:recursive}.
\end{proof}

\begin{prop} \label{prop3^6}
There exists an SBGDD of type $3^6$.
\end{prop}

\begin{proof}
We apply Construction~\ref{main-cons} by decomposing $K_6 \cdot \overline{K_3}$ into the direct product $K_6 \times K_3$ and three copies of $K_6$ on the independent sets of the product.
In \cite{DN}, it is shown that there exists a multiset of triangles in $K_6 \times K_3$ that cover the 
$|E(K_6 \times K_3)|=90$ edges with frequencies $0,1,2,\dots,89$.  To this, we include SBTS$_\mu(6)$, $\mu=90,105,120$ on each of the sets of size 6.  This produces the desired SBGDD.
\end{proof}

\begin{prop} \label{prop5^8}
There exists an SBGDD of type $5^8$.
\end{prop}

\begin{proof}
Take a $\{4,5\}$-GDD of type $5^8$ and replace blocks with SBTS$_\mu(4)$ or SBTS$_\mu(5)$.  We note that not all blocks have size four; hence, 
we may begin the interval of frequencies using at least one SBTS$_\mu(5)$, and finish using SBTS$_\mu(4)$ with $\mu >0$.
\end{proof}

\section{Combining the constructions}
\label{sec:combining}

We begin by highlighting a few instances of Construction~\ref{main-cons} carried by well known designs.  This sets the stage for the proof of Theorem~\ref{main} in two steps, first for $u \in \{4,5,6,8\}$, and then for other values of $u$.

\subsection{TD template}

The idea here is to use a transversal design as a template, or `master' design, replacing blocks with SBTS$_\mu$ (or minor variations thereof).

\begin{lemma}
\label{td}
Suppose there exists a TD$(k,n)$, where $n>1$ and $k \in \{4,5,6,8\}$.  Then there exists an SBGDD of type $n^k$.
\end{lemma}

\begin{proof}
We replace each block $B$ of the TD with an SBTS$_\mu(|B|)$ so that the frequency intervals abut.  In the case $k=4$, we take care to apply the adjustment outlined in Section~\ref{sec:recursive} to the first five blocks.  For $k \in \{5,8\}$, we apply the frequency interleaving technique to three blocks at a time, as detailed in Section~\ref{sec:recursive}.
\end{proof}

\rks
This result holds with the same proof for any $k \equiv 0,1 \pmod{3}$, but for our application to follow it is enough to focus on the indicated block sizes.  Note that extending the proof to $k \equiv 2 \pmod{3}$, $k \ge 11$, would require an interleaving strategy for larger blocks.

We also remark that incomplete transversal designs can be used in a minor variant of this construction for situations when a TD$(k,n)$ is unavailable.  For instance, an SBGDD of type $6^4$ can be constructed from an SBGDD of type $2^4$ and a pair of orthogonal latin squares missing a common $2 \times 2$ subsquare.  More details are given in \cite{Joey}.

\subsection{GDD inflation}

We have the following analogue of Wilson's fundamental construction.

\begin{lemma}
\label{inflate}
Suppose there exists a $\{3,4,5\}$-GDD of type $g^u$. Suppose further that, for each block size $k$ present in this GDD, there exists an SBGDD of type $m^k$ which has a sub-SBTS$(5)$ in the event that $k=5$ and $3 \nmid m$.  Then there exists an SBGDD of type $(mg)^u$.
\end{lemma}

\begin{proof}
Let $(X,\cB)$ be the hypothesized GDD.
Give every point weight $m$, so that the resultant set of points is $X \times [m]$. 
For each block $B \in \cB$, say with $|B|=k$, we replace it with an SBGDD of type $m^k$ on $B \times [m]$, where in the case $k=5$ and $3 \nmid m$ we use instead an SBIGDD of type $(m;1)^5$.

Note that we may lift frequencies of such designs using either a 3-GDD of type $m^{|B|}$ or 3-IGDD of type $(m;1)^5$, respectively.  In the case when $3 \nmid m$, we take care to first set aside and order the blocks of size $5$, replacing with SBTS$(5)$ and applying the interleaved lifting method in Section~\ref{sec:recursive} so that frequencies from such blocks form an interval.  Then, we ensure that the SBGDDs of type $m^3$, $m^4$, $m^5$ (or SBIGDD$_\mu(m;1)^5$) have abutting frequency lists.
\end{proof}

\rks
Note that we may take $m=1$ in the event that no blocks of size three occur in the GDD.  Or, if the GDD has exclusively blocks of size three, then the conclusion holds automatically for any $m \ge 2$ in view of the existence \cite{DN} of Sarvate-Beam cubes of order $m$.


\subsection{Proof of the main result}

As noted earlier, the case $g=1$ was settled in \cite{DSG,MCF} and the case $u=3$ was settled in \cite{DN}.  The remaining cases are handled by the following two results.

\begin{prop}
\label{exceptional-u}
There exists an SBGDD of type $g^u$ for any $g \ge 2$ and $u \in \{4,5,6,8\}$. 
\end{prop}

\begin{proof}
Suppose first that $u \in \{4,5\}$.  Lemma~\ref{td} produces an SBGDD of type $g^u$ except 
when the corresponding TD is not available, namely for $g \in \{2,6\}$ or in the case $10^5$. 
When $g=2$, we have a direct construction.  For $6^4$, we apply Lemma~\ref{inflate} to a 3-GDD of type $2^4$, giving weight $3$.  For $6^5$, we apply Lemma~\ref{inflate} to a 4-GDD of type $3^5$, giving weight 2.  Finally, for $10^5$, we apply Lemma~\ref{inflate} to a 5-GDD of type $5^5$ (i.e. an affine plane with one parallel class removed), giving weight 2.

Next, suppose $u \in \{6,8\}$.
We divide the argument according to the prime factorization of $g$.
Suppose first that every prime divisor of $g$ is at least $u-1$.
There exists a TD$(u,g)$ by Lemma~\ref{td-existence}(a), and hence an SBGDD of type $g^u$ by Lemma~\ref{td}.
On the other hand, suppose $g$ has a prime divisor $p<u-1$, so that
$(p,u) \in \{(2,6), (2,8), (3,6), (3,8), (5,8)\}$.
In each of these cases, we have constructed an SBGDD of type $p^u$.  Suppose now that $g=pm$, $m \ge 2$.
By Lemma~\ref{small-gdds}, there exists a $\{3,4,5\}$-GDD of type $p^u$.  Also, 
there exists an SBGDD of type $m^k$ for each $k \in \{3,4,5\}$ from \cite{DN} and the first part of the proof.
So it follows from Lemma~\ref{inflate} that there exists an SBGDD of type $g^u$.
\end{proof}

\begin{prop}
\label{prop:closure}
There exists an SBGDD of type $g^u$ for any $g \ge 2$ and $u =7$ or $u \ge 9$.
\end{prop}

\begin{proof}
For the stated values of $u$, there exists a PBD$(u,\{3,4,5\})$ (see \cite[IV.3.23]{handbook}) which is equivalently a $\{3,4,5\}$-GDD of type $1^u$.  Also, there exists an SBGDD  of type $g^k$ for each $k \in \{3,4,5\}$ by Proposition~\ref{exceptional-u}.  The result now follows by letting $g$ take the role of $m$ in Lemma~\ref{inflate}. 
\end{proof}

A summary of cases of the proof for $g>1$ is shown in Table~\ref{summary}.

\begin{table}[htbp]
\begin{tabular}{l|c|c|c|c|c|c}
$u \downarrow~~g \rightarrow$ & 2 & 3 & 5 & 6 & 10 & other \\
\hline
3~ & \multicolumn{6}{c}{Sarvate-Beam cubes} \\
\cline{2-7}
4& special &  $4$-GDD & TD &  GDD$\times 3$ & \multicolumn{2}{c}{TD} \\
\cline{2-7}
5~& special &  $4$-GDD & TD & GDD$\times 2$ & GDD$\times 2$  & TD \\
\cline{2-7}
6~& special &  special  & \multicolumn{4}{c}{TD or GDD$\times m$} \\
\cline{2-7}
8~& special &  $4$-GDD  &  $\{4,5\}$-GDD & \multicolumn{3}{c}{TD or GDD$\times m$} \\
\cline{2-7}
$7, \ge 9$& \multicolumn{6}{c}{from PBD$(u,\{3,4,5\}$} \\
\cline{2-7}
\end{tabular}
\caption{Summary of constructions for SBGDDs}
\label{summary}
\end{table}

\section{Larger starting frequencies}
\label{sec:mu}

Using 3-GDDs to raise frequencies, and with a few other minor variations to the ingredients, we can settle the existence of SBGDDs having any starting frequency $\mu>0$.  We first consider the recursive constructions, $u=7$ or $u \ge 9$, followed by special constructions for $u \in \{3,4,5,6,8\}$.
In general, whenever a recursive construction uses an SBTS$_0(k)$ to replace a block of size $k$, we may instead use an SBTS$_\nu(k)$; however, the latter object only exists \cite{DSG} when $3 \mid \nu$ or $k \equiv 0,1 \pmod{3}$.  

For starting frequencies $\mu \equiv 0 \pmod{3}$, the modification is straightforward. Replacing blocks of size $5$ in the cases $\mu \equiv  1,2 \pmod 3$ needs to be handled slightly differently.  We consider each congruence class in turn.

For $\mu \equiv 1 \pmod 3$, append one configuration to the first SBTS$(5)$, as in Table \ref{three5}, so that the pair frequencies covered are $1, 2, 3,\dots,9,$ and $12$.  To the second copy of SBTS$(5)$, similarly add nine copies of the configuration which covers frequencies $10,11, 13,14,\dots,20$. From there, continue to add blocks of size five in sets of three according to the main construction.  Adding $\mu -1$ copies of a $(5,3,3)$-design to each block of size five gives starting frequency $\mu$ for any $\mu \equiv 1 \pmod 3$.  This requires the remaining number of blocks to be a multiple of three.  If this is not the case, we have $u \equiv 2 \pmod 3$  and $\mu \equiv 1 \pmod 3$, so the necessary conditions imply that the SBGDD has $g \equiv 0 \pmod 3$. In this case, the blocks of size 5 are replaced after inflation by $g$, ensuring their frequency lists can be grouped three at a time.

For starting frequencies $\mu \equiv 2 \pmod 3,$ define a system $T= \{abd,	abe, acd, ace, bcd,	bce\}$ on the point set $\{a,b,c,d,e\}$, which has one pair $de$, of frequency $0$ and all other pairs of frequency $2$. Two copies of $T$ are used alongside the configuration in Table \ref{three5} to cover frequencies from $\mu$ to $\mu + 29$. The construction for the first three blocks is summarized in Table  \ref{three5higher} with $\mu =2$. For starting frequency $\mu = 2+3k$, a $(5,3,3k)$-design can accompany each block of size 5. 
Subsequent blocks of size 5 are handled in a similar fashion to those listed in the table, and lifted to cover the next interval of frequencies using copies of a $(5,3,30)$-design. Similar to the case with $\mu \equiv 1 \pmod 3$, this construction requires the number of blocks of size 5 to be divisible by three. If instead the number of blocks of size 5 is congruent to 1 (mod 3), then $u \equiv 2 \pmod 3$. As $\mu \equiv 2 \pmod 3$, the necessary conditions ensure that $g \equiv 0 \pmod 3$, and we may proceed as in the previous case.

 \begin{table}[htbp]
 $$\begin{array}{c|cccccccccc}
\text{pair} & ab & ac & ad & ae & bc & bd & be & cd & ce & de \\
\hline
\text{SBTS} & 0 & 1 & 3 & 2 & 5 & 4 & 7 & 6 & 8 & 9\\
\text{two configurations} & 2 & 2 & 2 & 2 & 2 & 2 & 2 & 2 & 4 & 4\\
\hline
\text{totals} & 2&	3&	5	&4&	7	&6	&9	&8&	12	& 13\\
\end{array}$$

$$\begin{array}{c|cccccccccc}
\text{pair} & ab & ac & ad & ae & bc & bd & be & cd & ce & de \\
\hline
\text{SBTS} & 0 & 1 & 3 & 2 & 5 & 4 & 7 & 6 & 8 & 9\\
\text{two configurations} & 2 & 2 & 2 & 2 & 2 & 2 & 2 & 2 & 4 & 4\\
\text{two } T\text{-configurations} &2	&2	&4	&4	&4	&4&	4	&4	&4	&4\\				
\text{$(5,3,6)$-design} &6&	6	&6&	6	&6&	6	&6&	6	&6&6\\
\hline
\text{totals} &10&	11	&15&	14&	17	&16&	19&	18&	22&	23\\
\end{array}$$

$$\begin{array}{c|cccccccccc}
\text{pair} & ab & ac & ad & ae & bc & bd & be & cd & ce & de \\
\hline
\text{SBTS} & 0 & 1 & 3 & 2 & 5 & 4 & 7 & 6 & 8 & 9\\
\text{two } T\text{-configurations} &2	&2	&4	&4	&4	&4&	4	&4	&4	&4\\				
\text{$(5,3,18)$-design} &18&	18	&18&	18	&18&	18	&18&	18	&18&18\\
\hline
\text{totals}& 20& 21 &	25&	24&	27&	26&	29&	28&	30&	31\\
\end{array}$$
 \caption{Lifting blocks of size 5 in designs with higher starting frequency} 
 \label{three5higher}
 \end{table}

We now consider the necessary special constructions with $u \in \{3,4,5,6,8\}$.  First, the incidence triples of a latin square of order $g$ produce a $3$-GDD of type $g^3$.  Taking $\lambda$ copies gives a $(3,\lambda)$-GDD of type $g^u$.  This, together with the main result of \cite{DN}, settles the case of SBGDD$_\mu$ of type $g^3$.  Similarly, starting frequencies for type $g^u$ can be increased by $\lam$ whenever Lemma~\ref{zhu} gives a $(3,\lambda)$-GDD of type $g^u$.
This settles existence for types $2^4, 2^5, 2^6,2^8,3^5,4^8, 6^4, 6^5,$ and $10^5$.  The only remaining exceptions to handle are SBGDD$_\mu$ of types $3^6, 3^8,$ and $5^8$. 

An SBGDD$_\mu$ of type $3^6$ can be constructed similarly as for $\mu=0$ in Proposition \ref{prop3^6}.  To lift frequencies, we first observe that a $K_3$-decomposition of $K_3 \times K_6$ is obtainable as a set of incidence triples of an idempotent latin square of order six with diagonal entries removed.  Include $\mu$ copies of this collection of blocks, and then SBTS$_\nu(6)$ beginning at $\nu=90+\mu$, $105+\mu$ and $120+\mu$.
An SBGDD$_\mu$ of type $3^8$ can be constructed as in Proposition~\ref{prop3^u} using ingredient SBTS$_\nu(4)$ whose frequencies are increased as needed. A similar modification to Proposition~\ref{prop5^8} produces an SBGDD$_\mu$ of type $5^8$. 

The above modifications to our constructions result in an extension to general starting frequency $\mu$, which we summarize below.

\begin{thm}
\label{arb-mu}
For positive integers $g,\mu,u$ with $u \ge 3$, there exists an SBGDD$_\mu$ of type $g^u$ if and only if $3 \mid g$ or $3 \mid \mu$ or $u \equiv 0,1 \pmod{3}$.
\end{thm}

\section*{Appendix}

Two SBGDDs found with the aid of a computer are given below. The blocks are listed without braces or commas for readability. The number preceding each block denotes its multiplicity.

\begin{itemize}
\item
SBGDD of type $2^6$ with group partition $\{\{a,a'\}, \{b,b'\}, \{c,c'\}, \{d,d\}',\{e,e'\},\{f,f'\}\}$ 

\footnotesize
\begin{center}
\begin{tabular}{cccccccccc}
$ace $&$ 2adf $&$ 5aef $&$ 3bcf $&$ 4bde $&$ 9bef $&$ 4cde $&$ 7cdf$\\
$24ab'd'$&$	19ab'e'$&$	20ac'e'$&$	21ac'f'$&$	21ad'f'$&$	7ba'c'$&$	8ba'd'$&$	9bc'e'$&$	bd'e'$&$	8bd'f'$\\
$10be'f'$&$	9ca'e'$&$	10ca'f'$&$	11cb'd'$&$	11cb'f'$&$	12cd'e'$&$	2cd'f'$&$	13da'b'$&$	13da'c'$&$	2db'c'$\\
$12db'e'$&$	2db'f'$&$	13dc'f'$&$	12de'f'$&$	16ea'd'$&$	15ea'f'$&$	17eb'c'$&$	15eb'f'$&$	17ec'd'$&$	18fa'b'$\\
$18fa'e'$&$	20fb'c'$&$	2fb'e'$&$	17fc'd'$&$	18fd'e'$&$	8a'b'c'$&$	10a'b'd'$&$	a'b'e'$&$	a'b'f'$&$	4a'c'd'$\\
$11a'c'e'$&$	a'c'f'$&$	11a'd'f'$&$	9a'e'f'$&$	11b'c'f'$&$	9b'd'e'$&$	12b'e'f'$&$	10c'd'e'$&$	11c'd'f'$&$	3d'e'f'$
\end{tabular}
\end{center}

\normalsize 

\item
SBGDD of type $2^8$ with group partition $\{\{a,a'\}, \{b,b'\}, \{c,c'\}, \{d,d\}',\{e,e'\},\{f,f'\},\{g,g'\},\{h,h'\}\}$

\footnotesize
\begin{center}
\begin{tabular}{cccccccccccccc}
$48abc$&$	ace$&$	4adf$&$	48adb'$&$	46aef$&$	13aeg'$&$	13aeh'$&$	3afg$&$	afg'$&$	4afh$&$	afh'$\\
$ab'c'$&$	24ab'd'$&$	19ab'e'$&$	ab'f'$&$	48ac'd'$&$	18ac'e'$&$	23ac'f'$&$	2ad'e'$&$	22ad'f'$&$	49ae'f'$&$	ae'g'$\\
$36agh$&$	4agh'$&$	ag'h$&$	28ag'h'$&$	bce$&$	2bcf$&$	2bde$&$	48bde'$&$	7bef$&$	48bea'$&$	47bfc'$\\
$5ba'c'$&$	8ba'd'$&$	2bc'd'$&$	7bc'e'$&$	bc'f'$&$	53bd'f'$&$	10be'f'$&$	bgh'$&$	4bg'h$&$	2bg'h'$&$	4cde$\\
$53cdf$&$	48cef'$&$	48ca'b'$&$	8ca'e'$&$	10ca'f'$&$	10cb'd'$&$	10cb'f'$&$	59cd'e'$&$	cd'f'$&$	8cgh$&$	3cgh'$\\
$cg'h$&$	11cg'h'$&$	48dec'$&$	3deg$&$	3deh$&$	12dfg'$&$	12dfh'$&$	11da'b'$&$	13da'c'$&$	47da'f'$&$	2db'c'$\\
$12db'e'$&$	2db'f'$&$	13dc'f'$&$	12de'f'$&$	13dgh$&$	dgh'$&$	4dg'h$&$	6dg'h'$&$	16ea'd'$&$	14ea'f'$&$	17eb'c'$\\
$47eb'd'$&$	15eb'f'$&$	17ec'd'$&$	21egh$&$	4egh'$&$	eg'h$&$	13eg'h'$&$	19fa'b'$&$	fa'c'$&$	48fa'd'$&$	19fa'e'$\\
$19fb'c'$&$	61fb'e'$&$	17fc'd'$&$	18fd'e'$&$	29fgh$&$	fgh'$&$	4fg'h$&$	21fg'h'$&$	8a'b'c'$&$	10a'b'd'$&$	a'b'f'$\\
$4a'c'd'$&$	52a'c'e'$&$	2a'c'f'$&$	2a'd'f'$&$	9a'e'f'$&$	8a'e'g'$&$	8a'e'h$&$	8a'f'g$&$	8a'f'h'$&$	32a'gh'$&$	34a'g'h$\\
$5a'g'h'$&$	b'c'd'$&$	52b'c'f'$&$	2b'd'e'$&$	11b'e'f'$&$	2b'gh'$&$	5b'g'h'$&$	9c'd'e'$&$	11c'd'f'$&$	8c'gh'$&$	10c'g'h$\\
$5c'g'h'$&$	8d'e'g$&$	8d'e'h'$&$	7d'f'g'$&$	7d'f'h$&$	10d'gh'$&$	9d'g'h$&$	5d'g'h'$&$	16e'gh'$&$	18e'g'h$&$	5e'g'h'$\\
$26f'gh'$&$	24f'g'h$&$	5f'g'h'$&																	
\end{tabular}
\end{center}
\end{itemize}
\normalsize


\end{document}